\documentclass[final,3p,times]{elsarticle}

\usepackage{epsfig}
\usepackage{graphicx}
\usepackage{graphics}
\usepackage{amsmath}
\usepackage{epstopdf}
\usepackage{amsfonts}
\usepackage{amssymb}
\usepackage{amsthm}
\usepackage{mathrsfs}
\usepackage{subfigure}
\usepackage{wasysym}
\usepackage{caption}
\usepackage{color}
\usepackage[normalem]{ulem}
\usepackage[mathscr]{euscript}

\newtheorem{theorem}{Theorem}[section]
\newtheorem{lem}[theorem]{Lemma}
\newtheorem{algo}[theorem]{Algorithm}
\newtheorem{prop}[theorem]{Proposition}
\newtheorem{defn}[theorem]{Definition}
\newtheorem{rem}[theorem]{Remark}
\newtheorem{cor}[theorem]{Corollary}

\numberwithin{equation}{section}

\newcommand{\ZZ}{\mathbb{Z}}
\newcommand{\RR}{\mathbb{R}}
\newcommand{\NN}{\mathbb{N}}

\def\bbeta{\mbox{\boldmath $\beta$}}
\def\balpha{\mbox{\boldmath $\alpha$}}

\def\bgamma{\mbox{\boldmath $\gamma$}}

\newcommand\bh{\mathbf{h}}

\newcommand\bp{\mathbf{p}}

\newcommand\bu{\mathbf{u}}
\newcommand\bP{\mathbf{P}}

\newcommand\cB{\mathcal{BL}}
\newcommand\cL{\mathcal{L}}

\newcommand\ie{{\it\thinspace i.e.}}

\def\1o2{{\hbox{$\frac{1}{2}$}}}
\def\1o4{{\hbox{$\frac{1}{4}$}}}
\def\3o4{{\hbox{$\frac{3}{4}$}}}
\def\fc{{\hbox{$\phi_0,\phi_1,\psi_0,\psi_1$}}}

\def\sod{{\hbox{$\frac{s}{h_1}$}}}

\def\tod{{\hbox{$\frac{t}{h_2}$}}}

\def\2od{{\hbox{$\frac{2}{h}$}}}
\def\uo2d{{\hbox{$\frac{1}{2h}$}}}

\def\do4{{\hbox{$\frac{h}{4}$}}}
\def\3do4{{\hbox{$\frac{3h}{4}$}}}
\def\do42{{\hbox{$\frac{h}{2}$}}}

\newcommand\cC{\mathcal{C}}

\newcommand{\barr}{\begin{array}}
\newcommand{\earr}{\end{array}}

\journal{Journal of Mathematical Analysis and Applications}

\begin{document}

\date{\OurVersion}
\begin{frontmatter}

\title{Convergence analysis of corner cutting algorithms refining points\\ and refining nets of functions}

\author[label1]{Costanza Conti\corref{cor1}}
\ead{costanza.conti@unifi.it}

\author[label2]{Nira Dyn}
\ead{niradyn@post.tau.ac.il}

\author[label3]{Lucia Romani}
\ead{lucia.romani@unimib.it}

\cortext[cor1]{Corresponding author.}

\address[label1]{Dipartimento di Ingegneria Industriale, Universit\`{a} di Firenze, Italy}

\address[label2]{School of Mathematical Sciences, Tel-Aviv University, Israel}

\address[label3]{Dipartimento di Matematica e Applicazioni, Universit\`{a} di Milano-Bicocca, Italy}

\begin{abstract}
In this paper we give an elementary proof of the convergence of corner cutting algorithms refining points,
in case the corner cutting weights are taken from the rather general class of weights considered by Gregory and Qu (1996).
We then use similar ideas, adapted to nets of functions, to prove the convergence of corner cutting algorithms refining nets of functions, in case the corner cutting weights are taken from a stricter class of weights than in the refinement of points.
\end{abstract}

\smallskip


\end{frontmatter}

\section{Introduction}\label{sec:intro}
Carl de Boor proved the convergence of corner cutting algorithms refining points in a very general case \cite{deBoor87}.
Convergence of a wide class of corner cutting algorithms was proved by Gregory and Qu \cite{Gregory}.
In this paper we give a simple proof of the convergence result in \cite{Gregory}. Our proof is based on the simple but crucial observation that the polygonal line $\cL(\bu^{[k+1]},\bp^{[k+1]})$ through the points of level $k+1$, $(u_i^{[k+1]},p_i^{[k+1]}),\ i\in \ZZ$, is a piecewise linear interpolant of $\cL(\bu^{[k]},\bp^{[k]})$.
Using an elementary error formula, we show that the sequence of polygonal  lines $\{\cL(\bu^{[k]},\bp^{[k]})\}_{k \geq 0}$ is a Cauchy sequence, in case the corner cutting weights are taken from the rather general class of weights ${\mathscr W}$ considered in \cite{Gregory} and satisfy an asymptotic condition on their size.
We then adapt the approximation idea to show the convergence of corner cutting algorithms for bivariate nets of functions, when the weights are taken from ${\mathscr W}$ but satisfy a stricter asymptotic condition.
Both convergence results for points and for nets of functions are proved under a condition on the initial data.
Besides the theoretical interest in these two nice convergence results, corner cutting algorithms for nets of functions (points) generate a variety of $C^0$ bivariate functions ($C^0$ curves) approximating the initial net (polygonal line), with the corner cutting weights acting as shape parameters.
In a future work we plan to study the smoothness of the limits in the case of nets, and to derive conditions on the corner cutting weights which guarantee $C^1$ limit functions.
This was investigated in the case of points in \cite{deBoor90} and in \cite{Gregory}.\\

\noindent
The structure of the paper is as follows. In Section \ref{sec:cc_points} we give our proof of the convergence of corner cutting algorithms refining points (polygonal lines). In Section \ref{sec:cc_nets} we consider the case of bivariate nets of functions. First, in Subsection \ref{subsec:cc_nets1} we give preliminary results on Coons patches (see \cite{C64}) and their approximation properties since they are analogous to linear interpolants in the case of points. Then, in Subsection \ref{subsec:cc_nets2} we introduce the notion of bivariate nets of functions and present the corner cutting algorithms for them. The convergence theorem and its proof are given in Subsection \ref{subsec:cc_nets3}.

\section{Corner cutting algorithms for points in $\RR^n$}\label{sec:cc_points}

Corner cutting algorithms for points are
iterative methods that starting from a given sequence of points  $\bp^{[0]}=\{p_i^{[0]}, i\in \ZZ\}$ produce at each iteration denser and denser sequences of points $\bp^{[k]}$, $k>0$.
Whenever convergent, 
they allow the user to define a continuous curve that approximates the shape described by the given polyline. Convergence of
corner cutting algorithms can be briefly defined as follows.

\begin{defn}\label{def:1}
A 
corner cutting algorithm is said to be convergent if, for any initial sequence $\bp^{[0]}=\{p_i^{[0]}, i\in \ZZ\}$, there exists
a function $F_{\bp^{[0]}} \in C(\RR)$ such that
$$
\lim_{k\rightarrow +\infty}\sup_{i\in \ZZ}|F_{\bp^{[0]}}(2^{-k}i)-p_i^{[k]}|=0.
$$
\end{defn}

\smallskip \noindent In this section we investigate the convergence of univariate corner cutting schemes assuming the corner cutting weights to satisfy the same assumptions as Gregory and Qu \cite{Gregory}.

\begin{defn}[Corner cutting weights]\label{def:alfabeta}
Let $\ell(\ZZ)$ be the set of scalar valued sequences indexed by $\ZZ$.
We denote by ${\mathscr W}$ a subset of $\ell(\ZZ) \times \ell(\ZZ)$ of the form
\begin{equation}\label{def:parameters}
{\mathscr W}:= \left \{ (\balpha, \bbeta) \in \ell(\ZZ) \times \ell(\ZZ) \, : \ \inf_{i \in \ZZ} \, \{ \alpha_i, \, 1-\beta_i, \, \beta_i-\alpha_i \}>0 \right \}.
\end{equation}
Moreover, for $\bgamma:=(\balpha, \bbeta) \in {\mathscr W}$ we define
\begin{equation}\label{def:mu_def}
\mu(\bgamma):=\sup_{i \in \ZZ} \, \{\beta_i-\alpha_i,\  1-\beta_{i-1}+\alpha_i\}.
\end{equation}
\end{defn}

\smallskip \noindent
Now let $\ell^n(\ZZ)$ denote the set of vector valued sequences indexed by $\ZZ$ and let $\bP=\{P_i\in \RR^n, \ i\in \ZZ\}\in \ell^n(\ZZ)$.
In the following we define the corner cutting operator for an arbitrary sequence $\bP$ of points in $\RR^n$.

\begin{defn}[Corner cutting operator]\label{def:CCoperator}
The \emph{corner cutting} operator with corner cutting weights $\bgamma:=(\balpha, \bbeta) \in {\mathscr W}$, denoted by $CC_{\bgamma}$, maps $\ell^n(\ZZ)$ into $\ell^n(\ZZ)$.
For $\bP\in \ell^n(\ZZ)$
\begin{equation}\label{eq:CC}
\begin{array}{c}
\\  \left(CC_{\bgamma}(\bP)\right)_{2i}=(1-\alpha_i)P_i+\alpha_iP_{i+1},\quad \quad  \left(CC_{\bgamma}(\bP)\right)_{2i+1}=(1-\beta_i)P_i+\beta_iP_{i+1}.
\end{array}
\end{equation}
\end{defn}

\begin{rem}
The corner cutting operator given in Definition \ref{def:CCoperator} is the same as the one studied in \cite{Gregory}.
A more general corner cutting operator is considered in \cite{deBoor87}. The condition required in \eqref{def:parameters}
on the corner cutting weights follows from the observation that
$$
\alpha_i=\frac{\| Q_{2i}-P_i \|_2}{\|P_{i+1}-P_i\|_2}, \qquad 1-\beta_i=\frac{\| P_{i+1}-Q_{2i+1} \|_2}{\|P_{i+1}-P_i\|_2}, \qquad \beta_i-\alpha_i=\frac{\| Q_{2i+1}-Q_{2i} \|_2}{\|P_{i+1}-P_i\|_2},
$$
where $Q_{2i}=\left(CC_{\bgamma}(\bP)\right)_{2i}$, $Q_{2i+1}=\left(CC_{\bgamma}(\bP)\right)_{2i+1}$.
\end{rem}

\begin{figure}[h!]
\centering
\includegraphics[width=8.5cm]{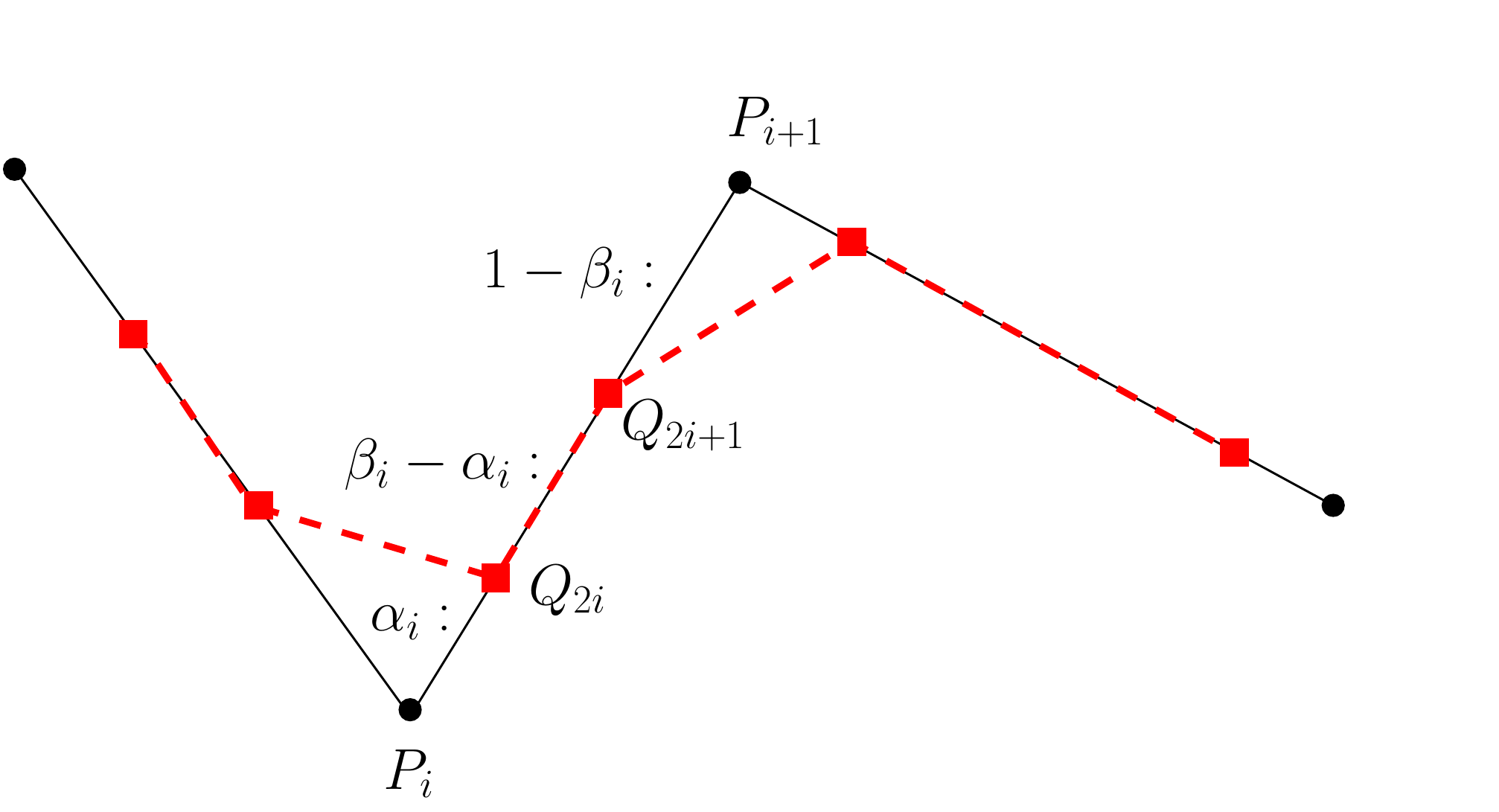}
\caption{One application of the $CC_{\bgamma}$-operator on a sequence of points in $\ell^2(\ZZ)$.
Here $Q_{2i}=\left(CC_{\bgamma}(\bP)\right)_{2i}$ and $Q_{2i+1}=\left(CC_{\bgamma}(\bP)\right)_{2i+1}$.}
\label{fig1}
\end{figure}

\smallskip \noindent
Denoting by $\bP^{[0]}\in \ell^n(\ZZ)$ a sequence of points in $\RR^n$ and assuming that, for each $k \ge 0$, a pair of scalar valued sequences $\bgamma^{[k]}:=(\balpha^{[k]},\bbeta^{[k]}) \in {\mathscr W}$ is assigned,
we can formulate the \emph{corner cutting algorithm}, for short the  $CC_{\bgamma}$-algorithm, as follows.

\begin{algo}
\label{def:CCalgorithm} Corner cutting algorithm for points:

{\sf \medskip \noindent  \ \ \ \ \ \ \ \ \ \ \ \ Input: $\bP^{[0]} \in \ell^n(\ZZ)$

\medskip \noindent  \ \ \ \ \ \ \ \ \ \ \ \ \ For $k=0,1,\dots,$

\medskip \noindent  \ \ \ \ \ \ \ \ \ \ \ \ \ \ \ \ \ \ \ Input:  $\bgamma^{[k]} \in {\mathscr W}$

\medskip \noindent  \ \ \ \ \ \ \ \ \ \ \ \ \ \ \ \ \ \ \ Compute \, $\bP^{[k+1]}=CC_{\bgamma^{[k]}}(\bP^{[k]})$} according to \eqref{eq:CC}
\end{algo}

\smallskip \noindent
In the remainder of this section we want to give a new simple proof of the fact that, for all choices of $\{\bgamma^{[k]} \in {\mathscr W}, \, k \ge 0\}$ satisfying $\sup_{k\ge 0} \, \mu(\bgamma^{[k]})<1$, and for all sequences of points in $\RR^n$ with bounded $L_{\infty}$ distance between every two consecutive points, the corner cutting algorithm always converges.
To this end we present two technical lemmas, where the first one is taken from \cite[Lemma 18]{ContiDyn11} and is here recalled for completeness.

\begin{lem}\label{lemma:18ContiDyn}
Let $f$ be a univariate function defined on $[a,b]$. If $f$ is Lipschitz continuous with Lipschitz constant $L$, then the error in approximating  $f$ by the linear interpolating polynomial at the points
$a,\ b$,
$$
{\cal L}(a,b;f(a),f(b))(x)=\frac{x-a}{b-a}f(b)+\frac{b-x}{b-a}f(a),$$
is bounded by
$$
\left| f(x)-{\cal L}(a,b;f(a),f(b))(x)\right|\le \frac{(b-a)L}2,\quad x\in [a,b]\,.
$$
\end{lem}

\proof It is well known that
\begin{equation}\label{eq:errornew}
f(x)-{\cal L}(a,b;f(a),f(b))(x)=(x-a)(x-b)[a,b,x]f
\end{equation}
with $[a,b,x]f$ the divided difference of order $2$ of $f$ at the points $a,b,x$. By definition of divided differences we get
\begin{equation}\label{eq:dd}
f(x)-{\cal L}(a,b;f(a),f(b))(x)= (x-a)(x-b)[a,b,x]f=\frac{(x-a)(x-b)}{b-a}\left(\frac{f(b)-f(x)}{b-x}-\frac{f(x)-f(a)}{x-a} \right)\,.
\end{equation}
Since $\frac{|(x-a)(x-b)|}{b-a} \le \frac{b-a}{4}$, and $f$ is Lipschitz continuous, then (\ref{eq:dd}) yields
$$
|f(x)-{\cal L}(a,b;f(a),f(b))(x)|\le \frac{(b-a)}{4}\left(\frac{L |b-x|}{|b-x|}+\frac{L|x-a|}{|x-a|}\right)=\frac{(b-a)L}2.
$$
\qed

\smallskip \noindent
The next lemma is about piecewise Lipschitz continuous functions.

\begin{lem}\label{Lem:PWL}
Let $f$ be a Lipschitz continuous function on each interval of a partition $\cdots <x_i<x_{i+1}<\cdots$ of the real line $\RR=\cup_{i\in \ZZ}[x_i,x_{i+1})$, with a bound $L$ on the Lipschitz constants.
Then $f$ is Lipschitz continuous in $\RR$ with Lipschitz constant $L$.
\end{lem}
\proof  Let  $t_1,\ t_2\in \RR$,  $t_1<t_2$. If $t_1,\ t_2$ belong to the same interval of the partition, say $t_1,\ t_2\in [x_i,x_{i+1})$, the inequality $|f(t_2)-f(t_1)| \le L \, |t_2-t_1|$ holds by assumption.
Otherwise, assuming $t_1\in [x_i,x_{i+1})$, $t_2\in [x_j,x_{j+1})$, $j\ge i+1$,  and writing
\begin{equation}
f(t_2)-f(t_1)=f(t_2)-f(x_j)+\sum_{l=i+1}^{j-1}\left(f(x_{l+1})-f(x_l)\right)+f(x_{i+1})-f(t_1)
\end{equation}
we easily arrive at
$$
|f(t_2)-f(t_1)|\le L\, |t_2-x_j|+L\,\sum_{l=i+1}^{j-1}|x_{l+1}-x_l|+L\, |x_{i+1}-t_1|=L\, |t_2-t_1|,
$$
which concludes the proof.
\qed

\begin{figure}[h!]
\centering
\includegraphics[width=13.5cm]{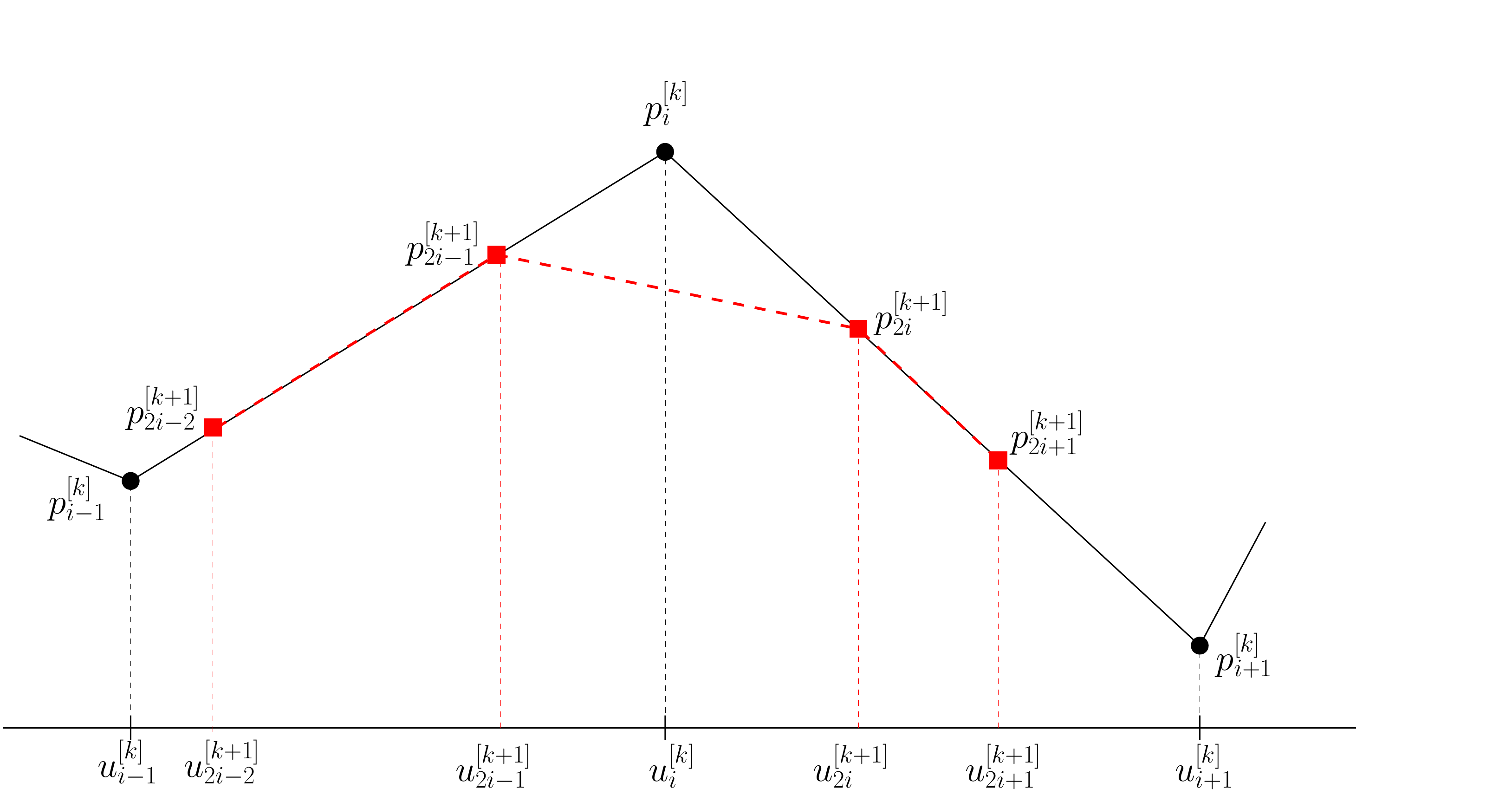}
\caption{$\cL(\bu^{[k+1]},\bp^{[k+1]})$ (dashed red) versus $\cL(\bu^{[k]},\bp^{[k]})$ (solid black).}
\label{fig2}
\end{figure}

\begin{theorem}\label{teo:CC-convergence}
For $\{\bgamma^{[k]}\}_{k\ge 0} \in {\mathscr W}$ such that
\begin{equation}\label{def:sup}
\sup_{k\ge 0} \, \mu(\bgamma^{[k]})<1,
\end{equation}
the corner cutting algorithm (Algorithm \ref{def:CCalgorithm}) converges for all initial sequences
$\bP^{[0]}=\{ P_i^{[0]} \in \RR^n, \, i \in \ZZ \}\in \ell^n(\ZZ)$ satisfying
$$
\|    P^{[0]}_{i+1}-P^{[0]}_i\|_{\infty}<L, \quad   \forall i\in \ZZ,
$$
with $L>0$.
\end{theorem}
\proof
We prove convergence of the $CC_{\bgamma}$-algorithm working component-wise.
First we introduce a parametrization at each refinement level.
Without loss of generality, we assume ${\bf u}^{[0]}=\ZZ$ and, for all $k \ge 0$, we denote by ${\bf u}^{[k]}$ the scalar sequence obtained from
${\bf u}^{[0]}$ by applying $k$ steps of the $CC_{\bgamma}$-algorithm (Algorithm \ref{def:CCalgorithm}).
Precisely, from the $(k-1)$-th level parameters, the $k$-th level parameters are obtained by the rules
$$
u_{2i}^{[k]}=(1-\alpha_i^{[k-1]})u_i^{[k-1]}+\alpha_i^{[k-1]} u_{i+1}^{[k-1]}, \qquad  u_{2i+1}^{[k]}=(1-\beta_i^{[k-1]})u_i^{[k-1]}+\beta_i^{[k-1]} u_{i+1}^{[k-1]}.
$$
Denoting by $p^{[k]}_i$ one component of $P^{[k]}_i$, we construct the piecewise linear interpolant to the data $(u_i^{[k]}, p_i^{[k]})$ and denote it by $\cL(\bu^{[k]},\bp^{[k]})$. In other words
$$
\cL(\bu^{[k]}, \bp^{[k]})(u)={\cal L}(u_i^{[k]},u_{i+1}^{[k]};p_i^{[k]},p_{i+1}^{[k]}),\quad u\in [u_i^{[k]}, u_{i+1}^{[k]}].
$$
By the assumption on $\bP^{[0]}$,  we know that  $|p^{[0]}_{i+1}-p^{[0]}_i|<L$ for all $i\in \ZZ$  and $\cL(\bu^{[0]},\bp^{[0]})$ is Lipschitz continuous with constant $L$ on $[u_i^{[0]}, u_{i+1}^{[0]}]=[i,i+1]$.
We show by induction that, for $k\ge 0$,  $\cL(\bu^{[k+1]},\bp^{[k+1]})$ is Lipschitz continuous with constant $L$ on $[u_i^{[k+1]}, u_{i+1}^{[k+1]}]$. Indeed, all points of $\bp^{[k+1]}$ lie
on $\cL(\bu^{[k]},\bp^{[k]})$ and therefore by the choice of ${\bf u}^{[k+1]}$  we know that
$|p_{i+1}^{[k+1]}-p_{i}^{[k+1]}| \le L\, |u_{i+1}^{[k+1]}-u_{i}^{[k+1]}|$.
Hence, by Lemma \ref{Lem:PWL}, we can conclude that $\cL(\bu^{[k]},\bp^{[k]})$ is Lipschitz continuous in $\RR$ with constant $L$ for all $k\ge 0$. Since  $\cL(\bu^{[k+1]}, \bp^{[k+1]})$ is, by construction, a piecewise linear interpolant to  $\cL(\bu^{[k]},\bp^{[k]})$, we can regard $\cL(\bu^{[k+1]},\bp^{[k+1]})$ as an approximation of $\cL(\bu^{[k]},\bp^{[k]})$. In particular, for
$u \in [ u_{2i}^{[k+1]}, u_{2i+1}^{[k+1]}]$, we have $|\cL(\bu^{[k+1]}, \bp^{[k+1]})(u)-\cL(\bu^{[k]},\bp^{[k]})(u)|=0$  (see Figure \ref{fig2}).
On the other hand, for $u\in [u_{2i-1}^{[k+1]}, u_{2i}^{[k+1]}],$
since $\cL(\bu^{[k]},\bp^{[k]})$ is Lipschitz continuous with constant $L$, we obtain by Lemma \ref{lemma:18ContiDyn} that
\begin{equation}\label{eq:approx_error}
|\cL(\bu^{[k+1]},\bp^{[k+1]})(u)-\cL(\bu^{[k]},\bp^{[k]})(u)|\le \frac12 \, L |u_{2i}^{[k+1]}-u_{2i-1}^{[k+1]}|\le \frac12 \, L \, d^{[k+1]},
\end{equation}
where $d^{[k]}=\sup_{i}|u_{i+1}^{[k]}-u_{i}^{[k]}|$.
Now, we proceed by comparing $d^{[k+1]}$ with $d^{[k]}$. To this purpose we have to distinguish between the following two cases (see Figure \ref{fig2}):

\begin{itemize}
\item Case 1: \ \ $u_{2i+1}^{[k+1]}-u_{2i}^{[k+1]}=(\alpha_{i}^{[k]}-\beta_{i}^{[k]})u_{i}^{[k]}+ (\beta_{i}^{[k]}-\alpha_{i}^{[k]}) u_{i+1}^{[k]}=(\beta_{i}^{[k]}-\alpha_{i}^{[k]}) (u_{i+1}^{[k]}-u_{i}^{[k]})$;

\item Case 2: \ \ $u_{2i}^{[k+1]}-u_{2i-1}^{[k+1]}=(1-\beta_{i-1}^{[k]})(u_i^{[k]}-u_{i-1}^{[k]})+\alpha_{i}^{[k]} (u_{i+1}^{[k]}-u_{i}^{[k]})$.
\end{itemize}
Both cases yield that $d^{[k+1]} \leq \mu^{[k]} \, d^{[k]}$ with $\mu^{[k]}:=\mu(\bgamma^{[k]})$.
Thus, in view of \eqref{eq:approx_error}, we get that
$|\cL(\bp^{[k+1]},\bu^{[k+1]})(u)-\cL(\bu^{[k]},\bp^{[k]})(u)| \leq \frac12 L d^{[k+1]} \leq \frac12 L \, d^{[0]} \, (\prod_{h=0}^{k} \mu^{[h]})$. Taking into account also that
$\prod_{h=0}^{k} \mu^{[h]} <\mu^{k+1}$ with $\mu:=\sup_{k\ge 0} \, \mu^{[k]}$, for any arbitrary $m \in \ZZ_+$ we can write
$$\begin{array}{ll}
|\cL(\bu^{[k+m]},\bp^{[k+m]})(u)-\cL(\bu^{[k]},\bp^{[k]})(u)| \leq &\sum_{\ell=0}^{m-1} |\cL(\bu^{[k+\ell+1]},\bp^{[k+\ell+1]})(u)-\cL(\bu^{[k+\ell]},\bp^{[k+\ell]})(u)| \\
\\
&\leq
\frac12 L \, d^{[0]} \mu^{k+1} \, \Big(\sum_{\ell=0}^{m-1} \mu^{\ell} \Big) \leq \frac{L d^{[0]}}{2(1-\mu)} \, \mu^{k+1},
\end{array}
$$
from which we conclude that
$\{\cL(\bu^{[k]},\bp^{[k]})\}_{k\ge 0}$ is a Cauchy sequence and therefore convergent. The limit of this sequence is the function $F_{\bp^{[0]}}$ of Definition \ref{def:1} (see, e.g., \cite{Dyn92}).
\qed

\begin{rem}
Some important observations:
\begin{itemize}

\item[(i)] The condition that the initial sequence of
points $\bP^{[0]}\in \ell^n(\ZZ)$ is such that $\|P^{[0]}_{i+1}-P^{[0]}_i\|_{\infty}<L$ for all $i\in \ZZ$, is equivalent to requiring
the piecewise linear interpolant to the data $(i,P_i^{[0]})$, $i\in \ZZ$  to be Lipschitz continuous with Lipschitz constant $L$.

\item[(ii)] Convergence of the corner cutting algorithm is not affected if in a finite number of iterations some or all of
the corner cutting weights $\alpha_i^{[k]}$ and $1-\beta_i^{[k]}$ are such that
$$ \alpha_i^{[k]}=0  \quad \hbox{or} \quad 1-\beta_i^{[k]}=0 \quad \hbox{or} \quad \beta_i^{[k]}-\alpha_i^{[k]}=0.$$
Even more,  convergence of the corner cutting algorithm is
guaranteed if in an infinite number of iterations 
the corner cutting weights satisfy \eqref{def:sup}.

\item[(iii)] Convergence of the corner cutting algorithm can be obtained under weaker assumptions on $\alpha_i^{[k]}$ and $\beta_i^{[k]}$
than the ones required in Theorem \ref{teo:CC-convergence}, namely by just requiring  that
$\sum_{\ell=0}^{m-1}\prod_{h=0}^{k+\ell} \mu^{[h]}<\infty$ for all $m \in \ZZ_+$.
\end{itemize}
\end{rem}

\section{Corner cutting algorithms for nets of functions}\label{sec:cc_nets}

The aim of this section is to show convergence of bivariate corner cutting algorithms refining nets of functions by suitably extending
the results introduced in the previous section.

\subsection{Preliminary results on the Coons pacth}\label{subsec:cc_nets1}

Since our proof of convergence of corner cutting schemes refining nets of univariate functions (u-functions for short) is based on error estimates for Coons interpolation, we need to recall first the definition of bilinear patch and  Coons patch (see \cite{G69}, \cite{Gordon82}). Then we point out some important properties of Coons patches that are relevant to our discussion.

\begin{defn}[The bilinear patch]\label{def:bilinear}
The bilinear patch interpolating the four points ${\cal P}=\{P_{ij}$, $i,j \in \{0,1\}\}$ is
$$ \cB({\cal P};\bh)(s,t)=
(1-\sod)\big((1-\tod)P_{00}+\tod
P_{01}\big)+\sod\big((1-\tod)P_{10}+\tod P_{11}\big),$$
where $\bh=(h_1,h_2)$ and $(s,t)\in [0,h_1]\times [0,h_2]$.
\end{defn}
\noindent
It is easy to verify that $$ \cB({\cal P};\bh)(i h_1,j h_2)=P_{ij}, \quad
i,j\in \{0,1\}.$$

\begin{defn}[The  Coons patch]\label{def:Coon}
Let $\phi_0(s), \phi_1(s)$, $s \in [0,h_1]$ and $\psi_0(t), \psi_1(t)$, $t \in [0,h_2]$  be four continuous univariate functions in $\RR^3$
such that $P_{ji}=\phi_i(j h_1)=\psi_j(i h_2)$ for $i,j \in \{0,1\}$.
The  Coons patch interpolating the four univariate functions $\phi_0, \phi_1$, $\psi_0, \psi_1$ is
\begin{equation}\label{Def:Coon}
\begin{array}{ll}
\cC(\fc;\bh)(s,t)=&(1-\sod) \psi_0(t)+ \sod \psi_1(t)+(1-\tod)
\phi_0(s)+ \tod \phi_1(s)- \cB({\cal P};\bh)(s,t),
\end{array}
\end{equation}
where $\bh=(h_1,h_2)$ and $(s,t)\in [0,h_1]\times [0,h_2]$.
\end{defn}

\smallskip \noindent
In the following, to simplify the notation we write $\cC(\phi,\psi;\bh)$ in place of $\cC(\fc;\bh)$.

\begin{rem}It is easy to verify the transfinite interpolation
properties of the Coons patch interpolant, \ie\ $$
\begin{array}{ll}
\cC(\phi,\psi;\bh)(0,t)=\psi_0(t),\quad \cC(\phi,\psi;\bh)(h_1,t)=\psi_1(t),\\
\\
\cC(\phi,\psi;\bh)(s,0)=\phi_0(s),\quad
\cC(\phi,\psi;\bh)(s,h_2)=\phi_1(s).\end{array}$$
\end{rem}

Next, the notion of mixed second divided difference of a bivariate function $F$ is introduced.
\begin{defn}\label{defnira_MSDD}
The mixed second divided difference (MSDD) of a bivariate function $F$ at the points $(\sigma_i,t_j)\in \RR^2$, $i,j\in \{1,2 \}$
is defined as
$$
[\sigma_1,\sigma_2;\tau_1,\tau_2]F=\frac{1}{(\sigma_1-\sigma_2)(\tau_1-\tau_2)}(F(\sigma_1,\tau_1)+F(\sigma_2,\tau_2)-F(\sigma_2,\tau_1)-F(\sigma_1,\tau_2)).
$$
\end{defn}

The following result expresses the error between a bivariate function $F$ and the Coons patch interpolating its boundary univariate functions.

\begin{prop}\label{pro:B}
Let $F$ be a bivariate continuous function defined on a rectangular domain $R=[a,b]\times[c,d]$, and denote by $\cC(F_{|\partial R})$ the  Coons patch interpolating $F_{|\partial R}$. Then
$$
F(s,t)-\cC(F_{|\partial R})(s,t)=
\frac{(s-a)(s-b)(t-c)(t-d)}{(b-a)(d-c)}([b,s;d,t]F-[s,a;d,t]F+[s,a;t,c]F-[b,s;t,c]F),\quad (s,t)\in R.
$$
\end{prop}
\proof
Let $(\cL_s(F))(s,t)=\frac{s-a}{b-a}F(b,t)+\frac{b-s}{b-a}F(a,t)$  and
$(\cL_t(F))(s,t)=\frac{t-c}{d-c}F(s,d)+\frac{d-t}{d-c}F(s,c)$. In view of \eqref{eq:errornew} we get
$$
((I-\cL_s)(F))(s,t)=\frac{(s-a)(s-b)}{b-a}\left(\frac{F(b,t)-F(s,t)}{b-s}-\frac{F(s,t)-F(a,t)}{s-a}\right),$$
and
$$
((I-\cL_t)(F))(s,t)=\frac{(t-c)(t-d)}{d-c}\left(\frac{F(s,d)-F(s,t)}{d-t}-\frac{F(s,t)-F(s,c)}{t-c}\right).$$
Moreover, since $\cC(F_{|\partial R})=\cL_s(F)+\cL_t(F)-\cL_t(\cL_s(F))$, we can also write 
$$
\begin{array}{ll}
F(s,t)-\cC(F_{|\partial R})(s,t)&=((I-\cL_t)(I-\cL_s)(F))(s,t)\\
\\
&=(I-\cL_t)\left(\frac{(s-a)(s-b)}{b-a}\left(\frac{F(b,t)-F(s,t)}{b-s}-\frac{F(s,t)-F(a,t)}{s-a}\right)\right).
\end{array}
$$
Therefore,
$$
\begin{array}{ll}
F(s,t)-\cC(F_{|\partial R})(s,t)&=\frac{(t-c)(t-d)}{d-c} \frac{(s-a)(s-b)}{b-a}
 \left(\frac{[b,s]F(\cdot,d)-[s,a]F(\cdot,d)}{d-t}-\frac{[b,s]F(\cdot,t)-[s,a]F(\cdot ,t)}{d-t}\right.\\
 \\
 &\left. -\frac{[b,s]F(\cdot,t)-[s,a]F(\cdot,t)}{t-c}+\frac{[b,s]F(\cdot,c)-[s,a]F(\cdot,c)}{t-c}\right)\\
 \\
 &=\frac{(t-c)(t-d)}{d-c} \frac{(s-a)(s-b)}{b-a}
 \left(\frac{[b,s](F(\cdot,d)-F(\cdot,t))}{d-t}-\frac{[s,a](F(\cdot,d)-F(\cdot,t))}{d-t}\right.\\
 \\
 &\left. -\frac{[b,s](F(\cdot,t)-F(\cdot,c))}{t-c}+\frac{[s,a](F(\cdot,t)-F(\cdot,c))}{t-c}\right)\\
 \\
&=\frac{(t-c)(t-d)}{d-c} \frac{(s-a)(s-b)}{b-a}\left([b,s;d,t]F-[s,a;d,t]F+[s,a;t,c]F-[b,s;t,c]F
\right).
\end{array}
$$
\qed


We introduce an important property of bivariate functions, which plays an important role in the convergence analysis of corner cutting schemes refining nets of functions.

\begin{defn}\label{defnira_BMSDD}
A bivariate function $F$ defined on $\Omega \subset \RR^2$ has the bounded MSDD property (BMSDD property) with constant $L$ in $\Omega$ if for any $\sigma_1,\sigma_2;\tau_1,\tau_2 \in \RR$ such that $(\sigma_i,t_j)\in \Omega$, $i,j \in \{1,2\}$,
satisfies
$$
|[\sigma_1,\sigma_2;\tau_1,\tau_2]F| \leq L.
$$

\end{defn}
Combining Definition \ref{defnira_BMSDD} with Proposition \ref{pro:B} we get

\begin{cor}\label{coro:14nira}
Let $F$ be a bivariate continuous function defined on a rectangular domain $R=[a,b]\times[c,d]$, and denote by $\cC(F_{|\partial R})$ the  Coons patch interpolating $F_{|\partial R}$.
If $F$ has the BMSDD property with constant $L$ in $R$, then
$$
\|F(s,t)-\cC(F_{|\partial R})(s,t)\|_{\infty} \leq L \frac{(b-a)(d-c)}{4},\quad (s,t)\in R.
$$
\end{cor}

\begin{rem} Note that Corollary \ref{coro:14nira} is a generalization of Lemma  \ref{lemma:18ContiDyn} to bivariate transfinite interpolation.
\end{rem}

\subsection{Corner cutting of nets of u-functions}\label{subsec:cc_nets2}
In this section we discuss a generalisation of the Chaikin type corner-cutting algorithm for nets of u-functions, that was presented in \cite{CD04} and \cite{CD10}. To this purpose we start by introducing the notion of net of u-functions.

\begin{defn}[Net of u-functions]\label{def:Compatible_Net}
A \emph{net} $N$ is a continuous bivariate function defined on a \emph{grid of lines}
\begin{equation}\label{eq:T}
T=T((\bh^{[s]},\bh^{[t]}),O)=\{s_i \times \RR, \, i \in \ZZ \} \, \cup \, \{\RR \times t_j, \, j \in \ZZ \}\,,
\end{equation}
with $\bh^{[s]}, \bh^{[t]}$ bi-infinite sequences of positive numbers, $O=(s_0,t_0)$, $s_{i+1}=s_i+h^{[s]}_i,\ i=0,1,2,\dots\ $,
$s_{i-1}=s_i-h^{[s]}_{i-1},\ i=0,-1,-2,\dots$,
and similarly for $\{t_j\}_{j\in \ZZ}$ with $h^{[t]}$ replacing $h^{[s]}$.
In other words, $N$ consists of the continuous univariate functions (u-functions) $N(s,t_j)$ and $N(s_j,t)$, $j\in \ZZ$ defined on $\RR$.
The point $O$ is termed the \emph{origin} of $T$ and the intervals $[s_j,s_{j+1}]$, $[t_j,t_{j+1}]$, $j \in \ZZ$ are termed \emph{grid intervals}.
\end{defn}

\medskip \noindent
To stress the relation between a net $N$ of u-functions and the corresponding grid of lines we use the notation $N=N(T)$.

\begin{defn}[$C^0$ net]
A net $N$ is termed a $C^0$ net if all the u-functions $\phi_j(s)=N(s,t_j)$, $\psi_j(t)=N(s_j,t)$, $j \in \ZZ$ are $C^0$ and compatible, \ie \, satisfy $\phi_i(s_j)=\psi_j(t_i)$ for all $i,j \in \ZZ$.
\end{defn}

\begin{defn}[Piecewise Coons patch]\label{def:PiecewiceCoon}
For a $C^0$ net $N$ consisting of the u-functions $\phi_j,\ \psi_j, j \in \ZZ$,
we denote by $\cC(N)$ the piecewise Coons patch interpolating it, which is locally
defined as
$$
\cC(N)(s,t)=\cC(\phi_i,\phi_{i+1},
\psi_j,\psi_{j+1}; h_{i,j})(s-s_i,t-t_j)\,, \quad (s,t) \in [s_i,s_{i+1}] \times [t_j, t_{j+1}],\quad i,j\in \ZZ
$$
with $h_{i,j}=(h_i^{[s]},h_j^{[t]})$, $h_i^{[s]}=s_{i+1}-s_{i}$, $h_j^{[t]}=t_{j+1}-t_{j}\ \ $   $i, j \in \ZZ$.
\end{defn}

%

\bigskip \noindent We remark that for a $C^0$ net $N(T)$, the net obtained by evaluating  $\cC(N)$ along the grid lines of any grid $\tilde T$, 
is also $C^0$, since $\cC(N)$ is continuous. Hence the following
iterative procedure is well defined.

\medskip
\begin{algo}\label{algo3} Corner cutting algorithm for nets of functions:

{\sf \medskip \noindent Input: \,  a $C^0$ net $N^{[0]}(T^{[0]})$ with $T^{[0]}=\ZZ^2$  (namely $s_{i}^{[0]}=t_{i}^{[0]}=i$, $i\in \ZZ$)

\medskip \noindent For \, $k=0,1,\dots$

\medskip \noindent \ \ \ \ \ \ \ Input:  $\bgamma^{[s],[k]}:=(\balpha^{[s],[k]}, \bbeta^{[s],[k]}) \in {\mathscr W}$ \quad and \quad
$\bgamma^{[t],[k]}:=(\balpha^{[t],[k]}, \bbeta^{[t],[k]}) \in {\mathscr W}$\\

\medskip \noindent \ \ \ \ \ \ \ Compute $s_{2i}^{[k+1]}=(1-\alpha_i^{[s],[k]})s_i^{[k]}+\alpha_i^{[s],[k]}s_{i+1}^{[k]}$ \quad and \quad $s_{2i+1}^{[k+1]}=(1-\beta_i^{[s],[k]})s_i^{[k]}+\beta_i^{[s],[k]}s_{i+1}^{[k]}$, for $\ i\in \ZZ$\\

\medskip \noindent \ \ \ \ \ \ \ Compute $t_{2j}^{[k+1]}=(1-\alpha_j^{[t],[k]})t_j^{[k]}+\alpha_j^{[t],[k]}t_{j+1}^{[k]}$ \quad and \quad $t_{2j+1}^{[k+1]}=(1-\beta_j^{[t],[k]})t_j^{[k]}+\beta_j^{[t],[k]}t_{j+1}^{[k]}$, for $\ j\in \ZZ$\\

\medskip \noindent \ \ \ \ \ \ \ Define $T^{[k+1]}=\{s_i^{[k+1]}\times \RR,\ i\in \ZZ\}\cup   \{\RR\times t_j^{[k+1]}, \ j\in \ZZ\}$\\

\medskip \noindent \ \ \ \ \ \ \ Compute $N^{[k+1]}= \cC(N^{[k]})|_{T^{[k+1]}}$ 
}
\end{algo}

\smallskip 
We denote the mapping from $N^{[k]}$ to $N^{[k+1]}$  in the above algorithm by $N^{[k+1]}= BC_{\bgamma^{[s], [k]},\bgamma^{[t],[k]}}(\cC(N^{[k]})$.
\noindent In the next subsection we prove that Algorithm \ref{algo3} is convergent under suitable assumptions on the initial net $N^{[0]}$ and the corner cutting weights.
To state the assumption on the initial net, we introduce the notion of a BMSDD net of functions which is a direct analogue  of Definition \ref{defnira_BMSDD}.

\begin{defn}A net of functions $N(T)$ has the BMSDD property with constant $L$, if
\begin{equation}\label{eq1nira}
| [\sigma_1,\sigma_2; \tau_1, \tau_2]N | \leq L,\quad \hbox{for}\quad  (\sigma_i,t_j) \in T,\quad i,j\in \{1,2\},
\end{equation}
where $$[\sigma_1,\sigma_2; \tau_1, \tau_2]N=\frac{1}{(\sigma_1-\sigma_2)(\tau_1-\tau_2)}(N(\sigma_1,\tau_1)+N(\sigma_2,\tau_2)-N(\sigma_2,\tau_1)-N(\sigma_1,\tau_2)).
$$
\end{defn}

\subsection{Convergence of the corner cutting algorithm for nets of functions}\label{subsec:cc_nets3}
In this subsection we state and prove the main result of this paper.

\begin{theorem}\label{teo:convergence}
Let $N^{[0]}$ be a $C^0$ net having the BMSDD property with constant $L$.
Then the corner cutting algorithm for nets of functions (Algorithm \ref{algo3}) is convergent for all $\{\bgamma^{[s],[k]}, \bgamma^{[t],[k]}\}_{k\ge 0} \in {\mathscr W}$ such that
\begin{equation}\label{eq:mu}
\mu^*=\sup_{k \geq 0} \, \max \{\mu(\bgamma^{[s],[k]}), \, \mu(\bgamma^{[t],[k]}) \}<\frac{\sqrt{3}}{3}.
\end{equation}
\end{theorem}

\noindent To prove this theorem we need several intermediate results. The first is an important observation about the BMSDD property of nets of functions.

\begin{lem}\label{lemma3nira}
Let $N(T)$ satisfy the inequality in \eqref{eq1nira} for
\begin{itemize}
\item[(a)] $t_j \leq \tau_1,\tau_2 \leq t_{j+1}$; $\sigma_1=s_i$, $\sigma_2=s_{i+1}$, $i,j \in \ZZ,$
\item[\hbox{or}]
\item[(b)] $s_i \leq \sigma_1,\sigma_2 \leq s_{i+1}$; $\tau_1=t_j$, $\tau_2=t_{j+1}$, $i,j \in \ZZ.$
\end{itemize}
Then $N(T)$ has the BMSDD property with constant $L$.
\end{lem}
\proof
Given $\sigma_1<\sigma_2$,  $\tau_1<\tau_2$ such that  $(\sigma_i,\tau_j) \in T$ for $i,j\in \{1,2\}$, there are two possibilities:
\begin{itemize}
\item[(i)] $\sigma_1=s_i$, $\sigma_2=s_{i+\ell}$, for some $i\in \ZZ$, $\ell\in \NN$ and $\tau_1,\tau_2\in \RR,$
\item[\hbox{or}]
\item[(ii)] $\tau_1=t_j$, $\tau_2=t_{j+\ell}$, for some $j\in \ZZ$, $\ell\in \NN$ and $\sigma_1,\sigma_2\in \RR.$
\end{itemize}
We consider case (i); the proof in case (ii) is similar. We prove that  the inequality in \eqref{eq1nira} holds for case (i) by induction. First we prove by induction on $\ell$ that  the inequality in \eqref{eq1nira}
holds in the case
$$
\hbox{(iii)}\qquad \quad  \sigma_1=s_i,\ \ \sigma_2=s_{i+\ell},\ \ \hbox{for some}\ i\in \ZZ,\ \ell\in \NN\ \  \hbox{and}\ \ t_j\le \tau_1,\tau_2\le t_{j+1}, \ \ \hbox{for some} \ j\in \ZZ.
$$
The above claim holds for $\ell=1$ by assumption ({\it a}). It remains to show that if  the inequality in \eqref{eq1nira} holds for $\ell\le m$, it holds for $\ell=m+1$. Now,
\begin{equation}\label{eqboh}
\begin{array}{ll}
[s_{i}, s_{i+m+1}; \tau_1, \tau_2]N  &=\frac{1}{(s_{i+m+1}-s_i)(\tau_2-\tau_1)} \left(N(s_{i+m+1},\tau_2)+ N(s_{i},\tau_1)-N(s_{i+m+1},\tau_1)-N(s_{i},\tau_2)  \right)\\
\\
&=\frac{(s_{i+m}-s_i)}{(s_{i+m+1}-s_i)}\frac{1}{(s_{i+m}-s_i)(\tau_2-\tau_1)} \left(N(s_{i+m},\tau_2)+ N(s_{i},\tau_1)-N(s_{i+m},\tau_1)-N(s_{i},\tau_2)  \right)\\
\\
&+\frac{(s_{i+m+1}-s_{i+m})}{(s_{i+m+1}-s_i)}\frac{1}{(s_{i+m+1}-s_{i+m})(\tau_2-\tau_1)} \left(N(s_{i+m+1},\tau_2)+ N(s_{i+m},\tau_1)-N(s_{i+m+1},\tau_1)-N(s_{i+m},\tau_2)  \right).
 \end{array}
\end{equation}
Thus by the induction hypothesis and by ({\it a}) we get for case (iii)
$$
|[s_{i}, s_{i+m+1}; \tau_1, \tau_2]N|\le \frac{(s_{i+m}-s_i)}{(s_{i+m+1}-s_i)}\, L +  \frac{(s_{i+m+1}-s_{i+m})}{(s_{i+m+1}-s_i)}\, L= L,
$$
and  the inequality in \eqref{eq1nira}  holds  in case (iii). This concludes the first part of the proof.

Next we prove, again by induction, that  the inequality in \eqref{eq1nira} holds in case (i).
We assume that  the inequality in \eqref{eq1nira} holds for $t_j\le \tau_1 \le t_{j+1}$ and $t_{j+m}\le \tau_2\ \le t_{j+m+1}$ for some $m\in \NN$, and show that
 the inequality in \eqref{eq1nira} holds for $t_j\le \tau_1 \le t_{j+1}$ and $t_{j+m+1}\le \tau_2 \le t_{j+m+2}$. This is sufficient since the case $m=0$ corresponds to case (iii).
Now, for $t_j\le \tau_1 \le t_{j+1}$ and $t_{j+m+1}\le \tau_2 \le t_{j+m+2}$
$$
\begin{array}{ll}
[s_{i}, s_{i+\ell}; \tau_1, \tau_2]N  &=\frac{1}{(s_{i+\ell}-s_i)(\tau_2-\tau_1)} \left(N(s_{i+\ell},\tau_2)+ N(s_{i},\tau_1)-N(s_{i+\ell},\tau_1)-N(s_{i},\tau_2)  \right)\\
\\
&=\frac{(t_{j+1}-\tau_1)}{(\tau_2-\tau_1)}\frac{1}{(s_{i+\ell}-s_i)(t_{j+1}-\tau_1)} \left(N(s_{i+\ell},t_{j+1})+ N(s_{i},\tau_1)-N(s_{i},t_{j+1})-N(s_{i+\ell},\tau_1)  \right)\\
\\
&+\frac{(\tau_2-t_{j+1})}{(\tau_2-\tau_1)}\frac{1}{(s_{i+\ell}-s_{i})(\tau_2-t_{j+1})} \left(N(s_{i+\ell},\tau_2)+ N(s_{i},t_{j+1})-N(s_{i+\ell},t_{j+1})-N(s_{i},\tau_2)  \right).
 \end{array}
$$
Thus,
$$
[s_{i}, s_{i+\ell}; \tau_1, \tau_2]N=\frac{(t_{j+1}-\tau_1)}{(\tau_2-\tau_1)} [s_{i}, s_{i+\ell}; t_{j+1}, \tau_1]N + \frac{(\tau_2-t_{j+1})}{(\tau_2-\tau_1)}[s_{i}, s_{i+\ell}; t_{j+1}, \tau_2]N.
$$
The MSDD in the first term above corresponds to case (iii), since $\tau_1,\ t_{j+1}\in [t_j,t_{j+1}]$, and the MSDD in the second term above corresponds to case (i) with $m$, since $t_{j+1}\in [t_{j+1},t_{j+2}]$ and $\tau_2\in [t_{j+1+m},t_{j+1+m+1}]$. By the first part of the proof we have $|[s_{i+\ell}, s_{i}; t_{j+1}, \tau_1]N|\le L$ and  by the induction hypothesis we have $[s_{i+\ell}, s_{i}; t_{j+1}, \tau_2]N\le L$. Thus,
$$
|[s_{i}, s_{i+\ell}; \tau_1, \tau_2]N|\le \frac{(t_{j+1}-\tau_1)}{(\tau_2-\tau_1)} L+ \frac{(\tau_2-t_{j+1})}{(\tau_2-\tau_1)}L= L,\quad \hbox{for}\quad t_j\le \tau_1 \le t_{j+1},\ \tau_{j+m+1}\le \tau_2 \le t_{j+m+2},
$$
and the inequality in \eqref{eq1nira} holds in case (i) for $m+1$.
\qed

\medskip \noindent A simple lemma follows from the linearity of the divided differences.

\begin{lem}\label{lemma4nira}
If $N_{\ell}(T)$, $\ell=1,...,m$ have the BMSDD property with constant $L$, then $\sum_{\ell=1}^m N_{\ell}(T)$ has the BMSDD property with constant $mL$.
\end{lem}

Using a similar induction to that in the second part of the proof of Lemma \ref{lemma3nira}, we can prove

\begin{lem}\label{lemma6nira}
A bivariate function which has the BMSDD property with constant $L$ on each rectangle of a grid $T$ has the BMSDD property with constant $L$ in $\RR^2$.
\end{lem}

The next lemma considers bivariate functions which are piecewise linear in one variable.

\begin{lem}\label{lemma7nira}
Let $F$ be a bivariate function of the form
$$F(s,t)=\frac{s-s_i}{s_{i+1}-s_i}F(s_{i+1},t)+ \frac{s_{i+1}-s}{s_{i+1}-s_i}F(s_i,t), \quad s \in [s_i,s_{i+1}],\ t \in \RR,\quad i \in \ZZ, $$
where $\{s_i\}_{i \in \ZZ} \subset \RR$ is an increasing sequence.
If $F$ satisfies  the inequality in \eqref{eq1nira} for $\sigma_1,\sigma_2 \in [s_i, s_{i+1}]$ for any $i \in \ZZ$ and $\tau_1, \tau_2 \in \RR$, then $F$ has the BMSDD property with constant $L$ in $\RR^2$.
\end{lem}

\proof
First we show that $F$ satisfies  the inequality in \eqref{eq1nira}  with constant $L$ in each rectangle of a grid $T$ defined by the parameters $\{s_i\}_{i \in \ZZ} $ and any increasing sequence $\{t_j\}_{j \in \ZZ}$.
Let $\sigma_1,\sigma_2 \in [s_i, s_{i+1}]$ and $\tau_1,\tau_2 \in [t_j, t_{j+1}]$ for some $i,j\in \ZZ$.
Since
$$[\sigma_1, \sigma_2; \tau_1, \tau_2]F=\frac{1}{(\tau_1- \tau_2)}\left(   [\sigma_1, \sigma_2]F(\cdot,\tau_1) -[\sigma_1, \sigma_2]F(\cdot,\tau_2)   \right),$$
the linearity of $F$ in $s$ implies that
$$
 [\sigma_1, \sigma_2]F(\cdot,\tau_j) =\frac{1}{(s_{i+1}- s_i)}\left(  F(s_{i+1},\tau_j) -F(s_i,\tau_j)   \right),\quad j=1,2,
 $$
and we get
$$[\sigma_1, \sigma_2; \tau_1, \tau_2]F=\frac{1}{(\tau_1- \tau_2)}\frac{1}{(s_{i+1}- s_i)}\left(  F(s_{i+1},\tau_1) +F(s_i,\tau_2)  - F(s_{i+1},\tau_2) -F(s_i,\tau_1)   \right).$$
From the assumption that $F$ satisfies the inequality in \eqref{eq1nira}  for $\sigma_1,\sigma_2 \in [s_i, s_{i+1}]$ for any $i \in \ZZ$ and $\tau_1, \tau_2 \in \RR$, we conclude that $F$ has the BMSDD
property with constant $L$ on each rectangle of $T$. Hence, by Lemma \ref{lemma6nira} $F$ has the BMSDD
property with constant $L$ in $\RR^2$.
\qed

\begin{rem}\label{rem8nira}
It is obvious that the same result holds if $F$ is linear in $t$ in each rectangle of $T$.
\end{rem}

Another important observation is

\begin{rem}\label{rem9nira}
The restriction to a grid of a bivariate function which has the BMSDD property with constant $L$ in $\RR^2$ is a net which has the BMSDD property with constant $L$.
\end{rem}

The next Theorem is our first key result.

\begin{theorem}\label{teo10nira}
If $N(T)$ has the BMSDD property with constant $L$ then $\cC(N)$ has the BMSDD property with constant $3L$.
\end{theorem}

\proof
Define the bivariate functions related to $N$  (similarly to the bivariate functions related to $F$ in the proof of Proposition \ref{pro:B})
$$\begin{array}{ll}
(\cL_s(N))(s,t)=&\frac{s-s_i}{s_{i+1}-s_i}N(s_{i+1},t)+\frac{s_{i+1}-s}{s_{i+1}-s_i}N(s_{i},t), \ s\in [s_i, s_{i+1}],\ t\in \RR, \ \ i\in \ZZ\\
\\
(\cL_t(N))(s,t)=&\frac{t-t_j}{t_{j+1}-t_j}N(s,t_{j+1})+\frac{t_{j+1}-t}{t_{j+1}-t_j}N(s,t_{j}), \ s\in \RR, \ t\in [t_j, t_{j+1}],\ \ j\in \ZZ,\\
\\
(\cL_s(\cL_t(N))(s,t)=&\frac{s-s_i}{s_{i+1}-s_i}\left(\frac{t-t_j}{t_{j+1}-t_j}N(s_{i+1},t_{j+1})+\frac{t_{j+1}-t}{t_{j+1}-t_j}N(s_{i+1},t_{j})\right)\\
&+\frac{s_{i+1}-s}{s_{i+1}-s_i}\left(\frac{t-t_j}{t_{j+1}-t_j}N(s_{i},t_{j+1})+\frac{t_{j+1}-t}{t_{j+1}-t_j}N(s_{i},t_{j})\right),\quad  \ s\in [s_i, s_{i+1}], \ t\in [t_j, t_{j+1}],\ \ i\in \ZZ,  j\in \ZZ .
\end{array}$$
Note that $\cL_s(\cL_t(N))$ is the piecewise bilinear function on the rectangles of $T$, interpolating the data $\{(s_i,t_j), N(s_i,t_j)\}_{i,j\in \ZZ}$.
It follows from \eqref{Def:Coon} and the definition of $\cC(N)$ that
\begin{equation}\label{eq7nira}
\cC(N)=\cL_s(N)+\cL_t(N)-\cL_s(\cL_t(N)).
\end{equation}
Next we show that the three functions in the right-hand side of the above equation have the BMSDD property with constant $L$ in $\RR^2$. By Lemma \ref{lemma7nira} and Remark  \ref{rem8nira} both
$\cL_s(N)$ and $\cL_t(N)$ have the BMSDD property with constant $L$ in each rectangle of $T$ since by assumption  $N(T)$  has the BMSDD property with constant $L$. Moreover, also $\cL_s(\cL_t(N))$
has the BMSDD property with constant $L$ in each rectangle of $T$ because for $\sigma_1,\sigma_2 \in [s_i, s_{i+1}]$ and $\tau_1,\tau_2 \in [t_j, t_{j+1}]$
$$[\sigma_1, \sigma_2; \tau_1, \tau_2]\cL_s(\cL_t(N))=[s_{i}, s_{i+1}; t_j, t_{j+1}]N.
$$
Now, by Lemma \ref{lemma6nira} the three functions have the BMSDD property with constant $L$ in $\RR^2$. Thus Lemma \ref{lemma4nira}, in view of \eqref{eq7nira}, implies that $\cC(N)$ has the BMSDD property with constant $3L$ in $\RR^2$.
\qed


A direct consequence of Theorem \ref{teo10nira} and Remark \ref{rem9nira} is

\begin{cor}\label{cor11nira}
If $N(T)$ has the BMSDD property with constant $L$, then $BC_{\bgamma^{[s]},\bgamma^{[t]}}(\cC(N))$ (defined after Algorithm \ref{algo3}) with $\bgamma^{[s]},\bgamma^{[t]} \in {\mathscr W}$, has the BMSDD property with constant $3L$.
\end{cor}

Corollary \ref{cor11nira} leads to our second key result.

\begin{cor}\label{cor12nira}
Let $\{N^{[k]}\}_{k \in \NN}$ be the nets generated by Algorithm \ref{algo3} from $N^{[0]}$.
If $N^{[0]}$ has the BMSDD property with constant $L$ then $N^{[k]}$ has the BMSDD property  with constant $3^kL$, for $k\ge 0$.
\end{cor}

\medskip We are now ready to prove the third key result.

\begin{theorem}\label{teo15nira}
In the notation of Algorithm \ref{algo3}, if $N^{[0]}$ has the BMSDD property with constant $L$ then
\begin{equation}\label{eq9nira}
\|\cC(N^{[k+1]})-\cC(N^{[k]})\|_{\infty} \leq 3^{k+1} L\frac{h_s^{[k+1]} h_t^{[k+1]} }{4} \end{equation}
where
$$
h_s^{[k+1]}=sup_{i \in \ZZ} (s_{i+1}^{[k+1]}-s_{i}^{[k+1]}), \qquad h_t^{[k+1]}=sup_{i \in \ZZ}(t_{i+1}^{[k+1]}-t_{i}^{[k+1]}).
$$
\end{theorem}

\proof
In view of Corollary \ref{cor12nira} and Theorem \ref{teo10nira}, $\cC(N^{[k]})$ has the BMSDD property with constant $3^{k+1}L$, and therefore by Remark \ref{rem9nira}, also $N^{[k+1]}=\cC(N^{[k]})|_{T^{[k+1]}}$ has this property. Regarding $\cC(N^{[k+1]})$ as the piecewise Coons patch interpolating $\cC(N^{[k]})|_{T^{[k+1]}}$, we conclude \eqref{eq9nira} from Corollary \ref{coro:14nira}.
\qed


\medskip We are now ready to prove Theorem \ref{teo:convergence}.

\medskip \noindent {\it Proof of Theorem \ref{teo:convergence}.}
By the way $T^{[k+1]}$ is constructed from $T^{[k]}$ in steps 1-3 of Algorithm \ref{algo3}, we see that
$$h_s^{[k+1]}\le \mu(\bgamma^{[s],[k]})h_s^{[k]}\quad \hbox{and}\quad h_t^{[k+1]}\le \mu(\bgamma^{[t],[k]})h_t^{[k]},$$
with $\mu(\bgamma^{[s],[k]})$ and  $\mu(\bgamma^{[t],[k]}) $ defined as in \eqref{def:mu_def}.
Defining $\mu^*=\sup_{k \geq 0} \, \max \{\mu(\bgamma^{[s],[k]}), \, \mu(\bgamma^{[t],[k]}) \}$ we get from \eqref{eq9nira}
$$
\|\cC(N^{[k+1]})-\cC(N^{[k]})\|_{\infty} \leq \frac{3^{k+1} L H}{4} (\mu^*)^{2k}=\frac{3L H}{4} (3(\mu^*)^2)^{k},$$
with $H=h^{[0]}_s\, h^{[0]}_t$.
Thus, if  $3(\mu^*)^2< 1$,  the sequence $\{\cC(N^{[k]}\}_{k \in \NN}$ is a Cauchy sequence and therefore convergent.
To conclude, the convergence of Algorithm \ref{algo3} is  guaranteed in case
 $$\mu^*=\sup_{k \geq 0} \, \max \{\mu(\bgamma^{[s],[k]}), \, \mu(\bgamma^{[t],[k]}) \}<\frac{\sqrt{3}}{3}.$$
 \qed

\begin{rem}The condition \eqref{eq:mu} in Theorem \ref{teo:convergence} can be relaxed to $$\sum_{k=0}^\infty 3^k \mu(\bgamma^{[s],[k]}), \, \mu(\bgamma^{[t],[k]}) <\infty.$$
\end{rem}

\bigskip
\bigskip
\noindent {\bf Acknowledgements.}
C. Conti and L. Romani acknowledge financial support from GNCS-INdAM.


\end{document}